\documentclass[12pt,oneside]{amsart}
\usepackage{amsmath}
\usepackage[mathscr]{euscript}
\usepackage{epsfig}
\usepackage{graphics}
\usepackage{amsbsy}
\usepackage{amsfonts}
\usepackage{amssymb}
\usepackage{graphicx}
\usepackage{rotating}
\usepackage{cancel}

 \newcommand{\bd}{\begin{definition}}
 \newcommand{\ed}{\end{definition}}
 \newcommand{\bt}{\begin{theorem}}
 \newcommand{\et}{\end{theorem}}
 \newcommand{\bp}{\begin{proposition}}
 \newcommand{\ep}{\end{proposition}}
 \newcommand{\bl}{\begin{lemma}}
 \newcommand{\el}{\end{lemma}}
 \newcommand{\bpr}{\begin{proof}}
 \newcommand{\epr}{\end{proof}}
 \newcommand{\bc}{\begin{corollary}}
 \newcommand{\ec}{\end{corollary}}

\usepackage{amssymb}
\usepackage{mathrsfs}
\newtheorem{theorem}{{\bf Theorem}}[section]
\newtheorem{lemma}[theorem]{{\bf Lemma}}
\newtheorem{corollary}[theorem]{{\bf Corollary}}
\newtheorem{proposition}[theorem]{{\bf Proposition}}

\newtheorem{definition}[theorem]{{\bf Definition}}
 \makeatletter
      \def\@setcopyright{}
      \def\serieslogo@{}
      \makeatother
\title{Basic Commutators in $n$-Lie Algebras}
\normalsize{
\author{Farshid Saeedi$^*$ and Seyedeh Nafiseh Akbarossadat}
}
\begin{document}
\maketitle

\noindent
\textbf{Abstract.} 
In this paper, we give the structure of free $n$-Lie algebras. Next, we introduce  basic commutators in $n$-Lie algebras and  generalize the Witt formula to calculate the number of the basic commutators. Also, we prove that the set of all of the basic commutators of weight $w$ and length $n+(w-2)(n-1)$ is a basis for $F^w$, where $F^w$ is the $w$th term of the lower central series in the free $n$-Lie algebra $F$.  
\ \\


\section{Introduction and History}

Lie polynomials appeared at the end of the $19${th} century and the beginning of the $20${th} century in the work of Campbell \cite{Campbell}, Baker \cite{Baker}, and Hausdorff \cite{Hausdroff} on exponential mappings in a Lie group, which has led to the so-called Campbell--Baker--Hausdorff formula. Around $1930$, Witt introduced the Lie algebra of Lie polynomials. He  showed that the Lie algebra of Lie polynomials is actually the free Lie algebra and that its enveloping algebra is the associative algebra of noncommutative polynomials. He proved what is now called the Poincaré--Birkhoff--Witt theorem and showed how a free Lie algebra is related lower central series of the free group. About at the same time,  Hall and Magnus, with their commutator calculus, opened the way to the bases of the free Lie algebra.  For more details about a historical review of free Lie algebras, we refer the reader to the reference \cite{[21]} and  its references.

The concept of basic commutators is defined in groups and Lie algebras, and there is also a way to construct and identify them. Moreover, a formula for calculating their number is obtained.

In 1962,  Shirshov \cite{Shirshov} gave a method that generalizes Hall's method for choosing a basis in a free Lie algebra. 

Basic commutators are of particular importance in calculating the dimensions of different spaces and are, therefore, highly regarded. Niroomand and  Parvizi \cite{Niroomand-Parvizi-M^2(L)} investigated some more results about the $2$-nilpotent multiplier $\mathcal{M}^{(2)}(L)$ of a finite-dimensional nilpotent Lie algebra $L$, and by using the Witt formula, they calculated its dimension. Moreover,  Salemkar,  Edalatzadeh, and  Araskhan  \cite{Salemkar-Edalatzadeh-Araskhan} introduced the concept of $c$-nilpotent multiplier $\mathcal{M}^{(c)}(L)$ of a finite-dimensional Lie algebra $L$ and obtained some bounds for $\mathcal{M}^{(c)}(L)$ by using the Witt formula and basic commutators. 

In 1985, Filippov \cite{vtf} introduced the concept of \textit{$n$-Lie algebras} as an $n$-ary multilinear  and skew-symmetric operation $[x_1,\ldots,x_n]$, which satisfies the following generalized Jacobi  identity:
\[[[x_1,\ldots,x_n],y_2,\ldots,y_n]=\sum_{i=1}^n[x_1,\ldots,[x_i,y_2,\ldots,y_n],\ldots,x_n].\]
Clearly, such an algebra becomes an ordinary Lie algebra when $n=2$.

In this paper, we are going to define  basic commutators in $n$-Lie algebras and give a formula to count  the number of them. 

\section{Free $n$-Lie algebras}
For the first time,  Akbarossadat and   Saeedi \cite{Akbarossadat-Saeedi-3} introduced the concept of free $n$-Lie algebras. In the following subsection, we will review the definition of this concept and some of its properties.

\subsubsection*{\underline{Free $n$-Lie Algebras}}

Given a set $X$, a free $n$-Lie algebra on $X$ over a field $\Bbb F$ is an $n$-Lie algebra $L$ over $\Bbb F$, together with a mapping $i:X\longrightarrow L$, with the following universal property: \\
For each $n$-Lie algebra, such as $M$, and each mapping $f:X\longrightarrow M$, there exists a unique $n$-Lie algebra homomorphism $F:L\longrightarrow M$ such that $f:F\circ i$. 

A standard argument shows that a free $n$-Lie algebra on $X$ is necessarily unique, up to $n$-Lie algebra isomorphism. Also, its existence is shown as follows:  

Let $X$ be a set. We construct a Lie algebra generated by $X$, satisfying no relations other than these:
\begin{eqnarray}
[x,x,\ldots, x]=0,\qquad \qquad \text{for all } x\in X,\label{eq1}\\
\sum_{i=1}^n[x_1,\ldots,x_{i-1},[x_i,x'_1,\ldots,x'_{n-1}],x_{i+1},\ldots,x_n]=0,\label{eq2}\\
\qquad\qquad\text{for all } x_i,x_j'\in X,\quad 1\leq i\leq n,~1\leq j\leq n-1.\nonumber
\end{eqnarray}
Let $M(X)$ be the set inductively defined as follows:
\begin{enumerate}
\item 
$ X\subseteq M(X)$.
\item 
If $x_1,\ldots,x_n \in M(X)$, then also $(x_1,\ldots,x_n)\in M(X)$.
\end{enumerate}
The set $M(X)$ is called the free magma on $X$. It consists of all bracketed expressions on the elements of $X$. On $M(X)$ we define an $n$-ary operation
\[\begin{array}{ccc}
 M(X)\times\cdots\times M(X)&\longrightarrow&M(X)\\ x_1\cdot~\ldots~\cdot x_n&\longmapsto&(x_1,\ldots,x_n).
\end{array}\]
For $m\subseteq M(X)$, we define its degree recursively: $deg(m) =1$ if $m \in X$ and $deg(m)=\sum_{i=1}^n m_i$ if $m=(m_1,\ldots,m_n)$. So the degree of an element $m\in M(X)$ is just the number of elements of $X$ that occur in $m$ (counted with multiplicities). For an integer $d\geq 1$, we set $M_d(X)$ as the subset of $M(X)$ consisting of all $m\in M(X)$ of degree $d$. Then 
\[M(X)=\bigcup_{d\geq 1}M_d(X).\] 
Let $\Bbb F$ be a field and let $A(X)$ be the vector space over $\Bbb F$ spanned by $M(X)$. If we extend the $n$-ary operation on $M(X)$ $n$-linearly to $A(X)$, then $A(X)$ becomes an (nonassociative) algebra; it is called the free $n$-algebra over $\Bbb F$ on $X$. Let $f\in A(X)$; if also $f\in M(X)$, then $f$ is said to be a monomial. For $f\in A(X)$, we define $deg(f)$ to be the maximum of the degrees of $m$, where $m$ runs over all $m\in M(X)$ that occur in $f$ with nonzero coefficient. 

Let $I_0$ be the ideal of $A(X)$ generated by all elements 
\begin{eqnarray*}
&(m, m,\ldots,m),\qquad \text{for all } m \in M(X),\\
&(m_1,\ldots,m_i,\ldots,m_j,\ldots,m_n) + (m_1,\ldots,m_j,\ldots,m_i,\ldots,m_n)\\&\qquad\qquad\qquad\qquad\qquad \text{for all } m_i\in M(X), ~ \text{for all } 1\leq i,j\leq n\\
&\sum_{i=1}^n(m_1,\ldots,m_{i-1},(m_i,m'_1,\ldots,m'_{n-1}),m_{i+1},\ldots,m_n),\qquad \text{for all } m_i,m'_i\in M(X).
\end{eqnarray*}
Set $L(X)=A(X)/I_o$. Let $\mathcal{B}$ be a basis of $L(X)$ consisting of (images of) elements of $M(X)$. Then it is immediate that we have 
\[(m, m,\ldots,m)=0,\]
\[(m_1,\ldots,m_i,\ldots,m_j,\ldots,m_n) + (m_1,\ldots,m_j,\ldots,m_i,\ldots,m_n)=0,\]
 and 
 \[\sum_{i=1}^n(m_1,\ldots,m_{i-1},(m_i,m'_1,\ldots,m'_{n-1}),m_{i+1},\ldots,m_n)=0,\]
  for all $m,m_i,m_i'\in\mathcal{B}$. It is easy to check that the relations  \eqref{eq1} and \eqref{eq2} hold for all elements of $L(X)$, so $L(X)$ is an $n$-Lie algebra. Therefore we will use a
bracket to denote the product in $L(X)$ and say $L(X)$ is the free Lie algebra on $X$. 

In other words, an $n$-Lie algebra $L$ is free on $X$ if and only if it is generated by all possible $n$-Lie bracketings of all elements of $X$. The only existing and possible relations among these bracketing generators are the  consequence of $n$-linearity bracketings, generalized Jacobi identity, and
\[[x_1, x_2, \dots, x_{i-1}, x_i, \dots, x_{j-1}, x_j,\dots, x_n] =-[x_1, x_2, \dots, x_{i-1}, x_j, \dots, x_{j-1}, x_i,\dots, x_n],\]
for all $1\leq i,j\leq n$.

The following lemma is a main tool for the proof of main theorem of this section.  
\bl \label{lemF^i,F^j}
Let $F$ be a free $n$-Lie algebra and let $i,j\in\Bbb N$. Then $[F^i,F^j,\underbrace{F,\dots,F}_{n-times}]\subseteq F^{i+j}$.  
\el
 \begin{proof}
We prove this by induction on $i$ and $j$. For each $i\geq 1$ and $j=1$, by the definition of lower central series terms, we have 
\[[F^i,F^j,F,\dots, F]=[F^i,F^1,F,\dots,F]=[F^i,F,F,\dots,F]=F^{i+1}.\]
Now, let $i$ be a fixed and arbitrary natural number and let the statement be true for $j-1$; that is, 
\[[F^i,F^{j-1},F,\dots, F]\subseteq F^{i+1}\qquad ({\text{induction hypothesis}}).\]
By the definition of lower central series terms $F^i$'s,  we know that $F^i=[F^{i-1},F,F,\dots,F]$. Hence every element of $F^i$ such an element $a$ is  $a=[[a_0,f_2',\dots,f_n'],f_2,\dots,f_n]$, where $a_0\in F^{i-1}$ and $f_s,f_r'\in F$, for each $2\leq r,s\leq n$. Now, let $x=[a,b,f_2,\dots,f_n]$ be an arbitrary element of $[F^i,F^j,\underbrace{F,\dots,F}_{n-times}]$, where $a=[a_0,f_2',\dots,f_n']\in F^i$, $b\in F^{j}$, and $f_s'\in F$, for all $2\leq s\leq n$. By the Jacobi identity, we have 
\begin{align*}
x=[a,b,f_2,\dots,f_n]&=[[a_0,f_2',\dots,f_n'],b,f_2,\dots,f_n]\\
&=[[a_0,b,f_2,\dots,f_n],f_2',\dots,f_n']+[a_0,[f_2',b,f_2,\dots,f_n],f_3',\dots,f_n']\\
&\quad+\dots +[a_0,f_2',\dots, f_{n-1}',[f_n',b,f_2,\dots,f_n]]\\
&=-\underbrace{[\underbrace{[b,a_0,f_2,\dots,f_n]}_{\in F^{j+i-1}},f_2',\dots,f_n']}_{\in F^{j+i}}-\underbrace{[\underbrace{[f_2',b,f_2,\dots,f_n]}_{\in F^{j+i-1}},a_0,f_3',\dots,f_n']}_{\in F^{j+i}}\\
&\quad+\dots +\underbrace{[a_0,f_2',\dots, f_{n-1}',\underbrace{[f_n',b,f_2,\dots,f_n]}_{\in F^{j+i-1}}]}_{\in F^{i+j}}\\ 
&=\sum_{s=2}^n[a_0,f_2',\dots,[f_s',b,f_2,\dots,f_n],f_{s+1}',\dots,f_n']\in F^{i+j}.
\end{align*} 
Thus $[F^i,F^j,\underbrace{F,\dots,F}_{n-times}]\subseteq F^{i+j}$.  
  \end{proof}

An ideal $I\unlhd L$ is called abelian, if $[I,I,L,\dots,L]=0$. Using this fact, we prove the following theorem, which is an important result. 
\bt \label{the1}
Let $F$ be a free $n$-Lie algebra and let $F^i$ be the $i$th term of the lower central series of $F$, for each $i\in\Bbb N$. Then $\dfrac{F^i}{F^{i+c}}$ is abelian, where $c=0,1,\ldots,i$.   
\et 
\begin{proof} 
It is obvious for $c=0$, since $\dfrac{F^i}{F^i}=0_{\frac{F^i}{F^i}}=F^i$. 
Let $c=1$. Since $[F^i,F^i,F,\dots,F]\subseteq [F^i,F,F,\dots,F]=F^{i+1}$, so $Z\left(\dfrac{F^i}{F^{i+1}}\right)=\dfrac{F^i}{F^{i+1}}$, and hence $\dfrac{F^i}{F^{i+1}}$ is abelian. 

Now suppose that $0\leq c\leq i$. We show that $\dfrac{F^i}{F^{i+c}}$ is abelian. In this regard, it is enough to show that $\dfrac{F^i}{F^{i+c}}$ is an abelian ideal of $\dfrac{F}{F^{i+c}}$, which is equivalent to 
\[[F^i,F^i,F,\dots,F]+F^{i+c}=0_{\frac{F^i}{F^{i+c}}}=F^{i+c},\qquad {\text{or}}\qquad [F^i,F^i,F,\dots,F]\subseteq F^{i+c}.\]
By Lemma \ref{lemF^i,F^j}, we know that $[F^i,F^i,F,\dots,F]\subseteq F^{2i}$. 
Since $0\leq c\leq i$, so $c+i\leq 2i$, and thus $F^{2i}\subseteq F^{i+c}$. Therefore, 
$[F^i,F^i,F,\dots,F]\subseteq F^{i+c}$. 
\end{proof}

\section{An introduction to basic commutators}
In this section, our aim is to define basic commutators of $n$-Lie algebras. For this, we need to introduce some concepts. 

\subsection*{\underline{Strings and collected parts}}
Let $X$ be a given set of noncommuting variables (also called letters).  We consider formal words or strings $b_1 b_2 \ldots b_n$, where each $b_i$ is one of the letters $x_1,x_2,\ldots,x_r,\ldots$. In other words, the word is a finite sequence of elements of $X$ or, equivalently, a noncommutative monomial. We include the empty word denoted $1$. We denote by $X^*$  the set of all words on $X$. With the product defined by the concatenation of words, $X^*$ becomes a monoid, which is the free monoid on $X$. Indeed, for any mapping $f:X\longrightarrow M$, where $M$ is any monoid, there is a unique homomorphism of monoids $F:X^*\longrightarrow M$ such that $f=F\circ i$, where $i$ is the natural injection $X\longrightarrow X^*$.  

A noncommutative polynomial is a linear combination of words with coefficients in  the field $\Bbb F$. The set of all these noncommutative polynomials becomes an algebra over $\Bbb F$ with the product inherited from the free monoid $X^*$, which is denoted by $\Bbb F\langle X\rangle$ and is the associative algebra $X$. Indeed, for any mapping $f:X\longrightarrow A$, where $A$ is any $\Bbb F$-algebra, there is a unique homomorphism of algebras as $F:\Bbb F\langle X\rangle\longrightarrow A$, such that $f=F\circ i$, where $i$ is the natural injection from $X$ into $\Bbb F\langle X\rangle$.   

Since $\Bbb F\langle X\rangle$ is an algebra, so from the product, we can define a bracket on it and turn it into a Lie algebra. For this, let $P$ and $Q$ be two  polynomials in $\Bbb F\langle X\rangle$. Define $[-,-]:\Bbb F\langle X\rangle\times \Bbb F\langle X\rangle\longrightarrow \Bbb F\langle X\rangle$ by $[P,Q]:=PQ-QP$. It is easy to check the bilinear and skew-symmetric properties and Jacobi identity. Hence this defines a structure of Lie algebra on $\Bbb F\langle X\rangle$.  

Now, let $\mathcal{B}_{\Bbb F}(X)$ denote the smallest Lie subalgebra of $\Bbb F\langle X\rangle$ that contains each letter in $X$. It is easy to prove that $\mathcal{B}_{\Bbb F}(X)$ is a free Lie algebra. In what follows, we show that every element of the Lie subalgebra $\mathcal{B}_{\Bbb F}(X)$ is a  Lie polynomial and that the set of all basic commutators on $X$ is a basis for $\mathcal{B}_{\Bbb F}(X)$. Hence the dimension of it is the number of all basic commutators on $X$. Indeed, for this purpose, we first need to introduce the concept of basic commutators and some of their properties. In what follows, we discuss this topic. 

Now, we introduce formal commutators $c_i$ and weights $\omega(c_i)$ by the rules: 

1) $c_i = x_i$, $i = 1, \ldots, r$, are the commutators of weight 1; that is, $\omega(x_i)=1$.

2) If $c_{i_1},\ldots, c_{i_n}$ are commutators, then $c_k=(c_{i_1},\ldots, c_{i_n})$ is a commutator and $\omega(c_k)=\sum_{j=1}^n\omega(c_{i_j})$.

Note that these definitions yield only a finite number of commutators of any given weight. We shall order the commutators by their subscripts, counting $c_i=x_i$, $i=1, \ldots, r$, and listing in order of weight, but giving an arbitrary ordering to commutators of the same weight.

A string $c_{i_1}\ldots c_{i_r}$ of commutators is said to be in collected form if  $i_1\leq\cdots\leq i_r$, that is, if the commutators are in order read from left to right. An arbitrary string of commutators $c_{i_1}\ldots c_{i_r}c_{i_{r+1}}\ldots  c_{i_s}$ will, in  general, has a collected part $c_{i_1}\ldots c_{i_r}$ if $i_1\leq\cdots\leq i_r$ and if $i_r\leq i_j$, $j=r+1, \ldots, s$. Moreover, this sting will have an uncollected part $c_{i_{r+l}} \ldots c_{i_s}$, where $i_{r+1}$ is not the least of $i_j$,  $j=r+1,\ldots,s$. 
The collected part of a string $c_{i_1}\ldots c_{i_s}$ will be void unless $i_1$ is the least of the subscripts.

\subsection*{\underline{Collecting process in strings}}

In this part, we define a collecting process for  strings of commutators. If $c_u$ is the earliest commutator in the uncollected part and if $c_{i_j}=c_u$ is the leftmost uncollected $c_u$, then we replace $c_{i_1}\ldots c_{i_m}\ldots c_{i_j-1}c_{i_j}\ldots c_{i_n}$ with $-c_{i_1}\ldots c_{i_m}\ldots c_{i_j}c_{i_j-1}\ldots c_{i_n}$. This has the effect of moving $c_{i_j}$ to the left as the formal coefficient of $(-1)$. Thus $c_{i_j}$ is still  the earliest commutator in the uncollected part. Alter  enough steps, $c_{i_j}$ will be moved to the $(m+1)$st position and will become a part of the collection part. Note that in every step of moving $c_{i_j}$ to the left, the coefficient of $(-1)$ is created that can be abridged.

\subsection*{\underline{Basic commutators and basic product}}

Suppose that sequences of basic commutators $c_1,c_2,\dots$ formed from the generators $x_1,x_2,\dots,x_d$ are given. We call a product of basic commutators $c_{i_1}c_{i_2}\dots c_{i_s}$, a {\it basic product} if it is in collected order; that is, $i_1\leq i_2\leq \dots\leq i_s$. The number of all basic products of weight $2$ and length $n$ formed from generators $x_1,x_2,\dots,x_d$, is $d^n$.

In applying the collecting process to a positive word, not all commutators will arise. Thus $(x_n,\dots,x_3,x_2,x_1)$ may arise but not $(x_n,\dots,x_1,x_2,x_3)$, since $x_3$ is collected before $x_1$. The commutators that may actually arise are called {\it basic commutators}. We give a formal definition of the basic commutators for $n$-Lie algebras generated by $x_1,\dots,x_d$.

\subsection*{\underline{Brief introduction of basic commutators in $n$-Lie algebras}}
In the previous sections, we completely and precisely define the basic commutators by using the words. Indeed in this section, we are going to briefly introduce the basic commutators for $n$-Lie algebra $L$ (especially) with the free representation $\frac{F}{R}$ over the field $\Bbb F$, in the bracket form. In other words, we introduce a method to determine the basic commutators in $n$-Lie algebras. Note that two properties of $n$-Lie algebras, as ``skew-symmetric" and ``Jacobian identity", are important and necessary to determine the basic commutators.

In general, it can be said that the set of all basic commutators $A$ is the smallest subset of the set of all commutators of $L$, named $W$, such that $A$ can be produced $W$.  

The basic commutators are characterized by two concepts as ``{\it weight}" and ``{\it length}". In the case when $n=2$, these two concepts are the same. The basic commutator $0\neq c\in F$ has weight $w_c$, if $c\in F^{w_c}$ and $c\not\in F^{w_c+1}$, where $F^{w_c}$ is $w_c$th terms of lower central series.
Also, the length of $c$ is the number of its components that are in $X$, and we denote it by $m_c$ and calculate it as follows:
\[m_c=n+(w_c-2)(n-1)=n+w_cn-2n-w_c+2=(w_c-1)n-(w_c-2).\] 
In the general case, weight and length are denoted by $w$ and $m$, respectively.  
Thus, the  general form of basic commutators with different weights and lengths is shown in Table \ref{tab1}.

\begin{table}[h!]
\caption{Generally form of basic commutators}\label{tab1}\begin{tabular}{cccc}\hline 
Domain& weight($w$)& length ($m$)& generally form\\ \hline 
$X$& 1&1& $x_i$\\ \hline 
$F^2$& 2& $n$& $[x_{i_1},\dots,x_{i_n}]$\\ \hline 
$F^3$& 3& $2n-1$& $[[x_{i_1},\dots,x_{i_n}],x_{j_2},\dots,x_{j_n}]$\\ \hline 
$\vdots$&$\vdots$&$\vdots$&$\vdots$\\ \hline
\end{tabular}
\end{table}

With the above discussion, we have  Table \ref{tab2}.

\begin{table}[h!]
\caption{The length of the basic commutators of weight $w$ in $n$-Lie algebras}\label{tab2}\begin{tabular}{cccccccccc}\hline
$w$&1&2&3&4&5&6&7&8&$\cdots$ \\ \hline
$n=2$&1 &2 &3 &4 &5 &6 &7 &8 & $\cdots$ \\ \hline
$n=3$&1 &3 &5 &7 &9 &11 &13 &15 & $\cdots$ \\ \hline
$n=4$&1 &4 &7 &10 &13 &16 &19 &22 & $\cdots$ \\ \hline
$n=5$&1 &5 &9 &13 &17 &21 &25 &29 & $\cdots$ \\ \hline
$n=6$&1 &6 &11 &16 &21 &26 &31 &36 & $\cdots$ \\ \hline
$n=7$&1 &7 &13 &19 &25 &31 &37 &43 & $\cdots$ \\ \hline
$n=8$&1 &8 &15 &22 &29 &36 &43 &50 & $\cdots$ \\ \hline
$\vdots$& $\vdots$& $\vdots$& $\vdots$& $\vdots$& $\vdots$& $\vdots$& $\vdots$& $\vdots$& $\vdots$ \\ \hline
\end{tabular}\end{table}
If in any $n$-Lie algebras (for each $n\in\Bbb N$), we denote the number of length of the basic commutator of weight $w_0$ by $m_{w_0}^n$, then we can obtain the following results:
\begin{itemize}
\item 
$m_{w_0}^n=m_n^{w_0}$,\qquad $n\geq 2$, and $w_0\geq 2$.
\item 
$m_{w_0}^3=m_{w_0}^2+m_{w_0-1}^2$.
\end{itemize}

Let $X=\{x_1,x_2,\ldots,x_n,\ldots\}$ be a basis set of $n$-Lie algebra $L$. We define the relation ''$<$'' on $X$ as $x_1<x_2<\cdots<x_n<\cdots$. In what follows, we define the basic commutators in $L$ as follows:
\begin{enumerate}
\item 
Every basic element $c_i^1=x_i$ in $X$ is the basic commutator of weight 1. 
\item 
The element $c_s^2=[x_{i_1},x_{i_2},\ldots,x_{i_n}]$ (where $x_{i_j}\in X$, for all $j=1,\ldots,n$) is a basic commutator of length $n$ and weight $2$. It belongs to  $F^2$, if $x_{i_1}>x_{i_2}>\cdots>x_{i_n}$; in other words, $i_1>i_2>\cdots>i_n$. We can compare every two basic commutators $c_s^2$ and $c_k^2$ with weight $2$ and length $n$. Suppose that $c_s^2=[x_{i_1},x_{i_2},\ldots,x_{i_n}]$ and that $ c_k^2=[x'_{j_1},x'_{j_2},\ldots,x'_{j_n}]$, and let $r_0$ be the first index (in moving from 
right to left) of $x_{i_r}$'s and $x'_{j_r}$'s which $x_{i_{r_0}}\neq x'_{j_{r_0}}$. If $x_{i_{r_0}}>x'_{j_{r_0}}$ or, equivalently, $i_{r_0}>j_{r_0}$, then we say  $c_{s}^2>c_k^2$, and otherwise, we say $c_{s}^2<c_k^2$. 
\item   
Let $c_t^w=[c_{j_1}^{w_1},\ldots,c_{j_n}^{w_n}]$, where all of $c_{j_1}^{w_1},\ldots,c_{j_n}^{w_n}$ are basic commutators of weights $w_1,\dots,w_n$ and lengths $n+(w_1-2)(n-1),\ldots,n+(w_n-2)(n-1)$, respectively. Then $c_t^w$ is called a basic commutator if the following conditions hold:
\begin{itemize}
\item 
$w>w_1\geq w_2\geq\cdots\geq w_n$.
\item 
Whenever $w_s=w_{s+1}$, then $c_{j_s}^{w_s}>c_{j_{s+1}}^{w_{s+1}}$, equivalently, $j_s>j_{s+1}$. 
\item 
Whenever $w_s>w_{s+1}$ and $c_{j_s}^{w_s}=[c_{j_{s_1}}^{w'_1},\ldots,c_{j_{s_n}}^{w'_n}]$, then $c_{j_{s_n}}^{w'_n}\leq  c_{j_n}^{w_n}<\cdots<c_{j_{s+2}}^{w_{s+2}}<c_{j_{s+1}}^{w_{s+1}}$. 
\end{itemize}
In this case, $c_t^w$ is an basic element of $F^w$ of weight $w$ and length $n+(w-2)(n-1)$. 
\end{enumerate}  
Note that every basic commutator of weight $w$ is smaller than every basic commutator of weight larger than $w$. Hence $c_i^1<c_j^2<c_k^3<\cdots$, for all indices $i,j,k,\ldots$ .
\subsection*{\underline{Collecting process in commutators of $n$-Lie algebras}}
It is obvious that there exist many commutators in generated free $n$-Lie algebra $F$ over $X$ such that they does not have a basic form. Indeed it is possible to choose an equivalent basic form for them  and characterize them by basic commutators. For this, we have to make the necessary changes to them by the properties of $n$-Lie algebras. 

Since $F$ is an $n$-Lie algebra, so every bracket with at least two same components is equal to zero. In fact, there is no nonzero bracket such that it has at least two same components in $F$, and all of the brackets in $F$ have $n$ different components.  

Suppose that $c=[c_{j_1}^{w_1},c_{j_2}^{w_2},\dots,c_{j_r}^{w_r},c_{j_{r+1}}^{w_{r+1}},\dots,c_{j_n}^{w_n}]$ is a nonbasic commutator in $F$. Then there are the following two cases: 

\begin{enumerate}
\item[(A)]
Let $c_{j_r}^{w_r}<c_{j_{r+1}}^{w_{r+1}}$. Then it is enough to write it as follows:
\[c=[c_{j_1}^{w_1},c_{j_2}^{w_2},\dots,c_{j_r}^{w_r},c_{j_{r+1}}^{w_{r+1}},\dots,c_{j_n}^{w_n}]=-[c_{j_1}^{w_1},c_{j_2}^{w_2},\dots,c_{j_{r+1}}^{w_{r+1}},c_{j_r}^{w_r},\dots,c_{j_n}^{w_n}].\] 
\item[(B)] 
Let $c_{j_1}^{w_1}=[{c'}_{r_1}^{w'_1},{c'}_{r_2}^{w'_2},\dots,{c'}_{r_n}^{w'_n}]$ be in the basic form. Then $c$ reads as follows:
\[c=[c_{j_1}^{w_1},c_{j_2}^{w_2},\dots,c_{j_n}^{w_n}]=[[{c'}_{r_1}^{w'_1},{c'}_{r_2}^{w'_2},\dots,{c'}_{r_n}^{w'_n}],c_{j_2}^{w_2},\dots,c_{j_n}^{w_n}].\]
We know that $c$ is a basic commutator if 
\[{c'}_{r_n}^{w'_n}\leq c_{j_n}^{w_n}< c_{j_{n-1}}^{w_{n-1}}<\dots< c_{j_3}^{w_3}< c_{j_2}^{w_2}.\]
Otherwise, $c$ is not in a basic form. Thus we replace $c$ by using the Jacobi identity and case (A). 
\end{enumerate}

Now, we state the theorems that are among the main results and our goals in this paper.

\bt \label{linear independent}
Let $F$ be a free $n$-Lie algebra with the ordered basis set $X=\{x_1,x_2,x_3,\dots\}$, where $x_1<x_2<x_3<\dots$. The set of all defined basic commutators of weight $w$ (denoted by $B(w)$) on $X$ is linearly independent. 
\et 
\begin{proof}  
It is obvious by the collecting process in commutators of $n$-Lie algebras. Since $F$ is the generated free $n$-Lie algebra on $X$, then it has no other relations except $n$-linearity, skew-symmetric properties, and Jacobi identity.
Therefore, all of the brackets in $F$ apply to these three properties. On the other hand, in the process of determining and identifying the basic commutators, we remove the commutators that are linearly dependent under these three properties ($n$-linearity, skew-symmetric properties, and Jacobi identity), and we separate only its linearly independent elements. Hence, all of the basic commutators of arbitrary and fixed weight $w$ are linearly independent.   
\end{proof}

The following theorem is the main result of this section, which is similarly discussed in group theory and Lie algebras.
\bt \label{maintheo1}
Let $F$ be a free $n$-Lie algebra. Then $\dfrac{F^i}{F^{i+c}}$ is abelian $n$-Lie algebras whose basis is the set of all basic commutators of weights $i, i+1, i+2, i+3,\dots, i+c-1$ and lengths $n+(i-2)(n-1),n+(i-1)(n-1),\ldots,n+(i+c-3)(n-1)$. 
\et 
\bpr
This theorem is proved similar to the Lie algebra case. 
\epr 

\section{Counting  of basic commutators}
It is necessary to note that counting  of basic commutators on every arbitrary basis set $X$ is not possible, in general. In fact, although it is possible to define the concept of basic commutators for any countable set, counting their number is only possible if the set $X$ is finite.  
The following formula, known as the {\it Witt formula}, is provided in group theory and Lie algebras for  counting of all of the basic commutators: 
\begin{equation}\label{Witt-formula}
l_d(w)=\dfrac{1}{w}\sum_{r|w}\mu(r)d^{\frac{w}{r}},
\end{equation}  
where $d$ is the number of generators of given group/Lie algebra (the number of members of $X$), $\mu$ is the {\it Möbius function}, which is defined for positive integers by the rule $\mu(1)=1$, and for $w=p_1^{a_1}p_2^{a_2}\dots p_s^{a_s}$ that $p_1,p_2,\dots, p_s$ are distinct primes numbers, by $\mu(p_1p_2\dots p_s)=(-1)^s$ and $\mu(p_1^{a_1}p_2^{a_2}\dots p_s^{a_s})=0$, if any $a_i>1$.

In group theory, the strategy of the proof of the above formula is to remove nonbasic commutators from all the commutators. That is, it first counts all possible commutators on $X$ and then removes nonbasic members and their periods. Indeed in $n$-Lie algebras, with this method, only a bound can be found for the number of basic commutators.

Now, let $X=\{x_1,x_2,\dots,x_d\}$ with the ordered relation $x_d>x_{d-1}>\dots>x_2>x_1$. We are going to determine the number of all of basic commutators of weight $w$ and length $m=n+(w-2)(n-1)$. Now, consider the following string of circles:  
\[\underbrace{\bigcirc\bigcirc\bigcirc\dots\bigcirc\bigcirc\bigcirc}_{m-times}\]
We want to put $x_i$'s in each circle such that the final result is a basic commutator. First, we consider all of the cases for putting them. It is easy to see that the number of these cases is $d^m$ (every circle has $d$ choices). Obviously, all of them are not basic, and hence we remove nonbasic cases. For this, we remove duplicate and linear dependent commutators by using two properties, skew-symmetric and Jacobi identity of $n$-Lie algebras. This method is similar to the method used in the group theory (see \cite{Hall}), and also, Shirshov \cite{Shirshov} provided it for Lie algebras. We know that
\[m=n+(w-2)(n-1)=n+wn-2n-w+2=(w-1)n-(w-2),\] 
and so $n$ is a divisor of $m$ if $w\stackrel{n}\equiv 2$; in other words, $w=qn+2$, for some integer positive $q$. Thus 
\[m=(w-1)n-(qn+2-2)=(w-1)n+qn=(w+q-1)n,\]
and hence $\dfrac{m}{n}=w+q-1$ is integer. Therefore, in general, $n\not| m$, but $n-1|m-1$, because 
\[m=n+(w-2)(n-1)=n+wn-2n-w+2,\]
and so 
\begin{align*}
m-1=n+wn-2n-w+2-1
&=\underline{n-2n+1}+\underline{wn-w}\\
&=-(n-1)+w(n-1)\\
&=(w-1)(n-1).
\end{align*}
We call a word $a_1a_2\dots a_s$ circular if $a_1$ is regarded as $a_s$, and also,  $a_1a_2\dots a_s$, $a_2a_3\dots a_sa_1$, $a_3\dots a_{n-1}a_na_2$, $\dots$, $a_na_1a_2\dots a_{n-1}$ are regarded as the same word. A circular word $C$ of weight $w$ and length $m=n+(w-2)(n-1)$ may conceivably be given by repeating a segment of $q$ letters $m/q$ times, where $q$ is some divisor of $m$. We say that $C$ is of period $q$ if we have 
\[\underbrace{(\underbrace{\bigcirc\bigcirc\dots\bigcirc}_{q-times})(\underbrace{\bigcirc\bigcirc\dots\bigcirc}_{q-times})\dots (\underbrace{\bigcirc\bigcirc\dots\bigcirc}_{q-times})}_{m-{\mbox{circles and}}~m/q-{\mbox{parts}}}.\]
So by the above, circular commutators must be removed from all of the generated commutators on $X$. Also, in weights that smaller than $w$, we have 
\[\underbrace{\underbrace{\bigotimes\bigotimes\bigotimes\dots\bigotimes}_{i-times}(\underbrace{\bigcirc\bigcirc\dots\bigcirc}_{q-times})(\underbrace{\bigcirc\bigcirc\dots\bigcirc}_{q-times})\dots (\underbrace{\bigcirc\bigcirc\dots\bigcirc}_{q-times})}_{(m-i)-{\mbox{circles and}}~(m-i)/q-{\mbox{parts}}}\]  
\begin{align}\label{generalized Witt-formula}
l_d^n(w)&\leq \dfrac{1}{m}\left(\sum_{r|m}\mu(r)d^{\frac{m}{r}}\right).
\end{align}  
Note that the equality holds if and only if $n=2$. 

In the next section, we introduce the recursive formula. 

\section[New method of counting of basic commutators]{New method of counting of basic commutators}
In the previous section, we provide the generalized Witt formula for counting basic commutators in $n$-Lie algebras with finite basis set $X$. Indeed, it could not discuss the relationship between the numbers of basic commutators of different weights. Now, in this section, we introduce some formulas in special cases with the new method and show their relationship. 

First, note the following theorem in which we give the number of all of the basic commutators of weight 2, in general $n$-Lie algebras with $d$ generators.
\begin{theorem}\label{weight2}
Let $L$ be an $n$-Lie algebra of dimension $d$. Then 
\[l_d^n(2)={{d}\choose{n}}.\]
\end{theorem}
\bpr
Let $L$ be an $n$-Lie algebra of dimension $d$ with the basis $X=\{x_1,x_2,\dots,x_d\}$ and put $w=2$. We are going to count the number of all basic commutators of weight $2$ on $d$ letters. It is clear that the number of all commutators of weight $2$ is equal to $d^n=n^n$. By the skew-symmetric property in $n$-Lie algebras, we have 
\[[x_{i_1},x_{i_2},\dots,x_{i_j},\dots,x_{i_k},\dots,x_{i_n}]=-[x_{i_1},x_{i_2},\dots,x_{i_k},\dots,x_{i_j},\dots,x_{i_n}],\] 
and hence $[x_{i_1},x_{i_2},\dots,x_{i_{j-1}},x,x_{i_{j+1}},\dots,x_{i_{k-1}},x,x_{i_{k+1}},\dots,x_{i_n}]=0$. Thus just independent and nonzero commutators of weight $2$ are as $[x_{i_1},x_{i_2},\dots,x_{i_n}]$, where all of $x_{i_j}$'s are distinct. So $l_d^n(2)={{d}\choose{n}}$. \epr 

\noindent{\bf Remark.} Since every $n$-Lie algebra has the skew-symmetric property, it is not possible to exist at least  two same components. Hence the number of basic commutators of weights greater than or equal to $2$ is equal to zero, if $d<n$.   
\bc \label{Cor-l_n^n(w)}
If $n=d$, then the following properties hold:
\begin{enumerate}
\item 
$l_n^n(1)=n$. 
\item 
$l_n^n(2)={{n}\choose{n}}=1$.
\item 
$l_n^n(3)={{n}\choose{n-1}}=n$. 
\item 
$l_n^n(4)
={{n}\choose{2}}+n$.
\item 
$l_n^n(5)
={{n}\choose{3}}+2{{n}\choose{2}}+n$.
\item 
$l_n^n(6)
={{n}\choose{4}}+3{{n}\choose{3}}+3{{n}\choose{2}}+n$.
\item 
$l_n^n(7)
={{n}\choose{5}}+4{{n}\choose{4}}+6{{n}\choose{3}}+4{{n}\choose{2}}+n$.
\item 
$l_n^n(8)
={{n}\choose{6}}+5{{n}\choose{5}}+10{{n}\choose{4}}+10{{n}\choose{3}}+5{{n}\choose{2}}+n$.
\item 
$l_n^n(9)
={{n}\choose{7}}+6{{n}\choose{6}}+15{{n}\choose{5}}+20{{n}\choose{4}}+15{{n}\choose{3}}+6{{n}\choose{2}}+n$.
\item 
$l_n^n(10)
={{n}\choose{8}}+7{{n}\choose{7}}+21{{n}\choose{6}}+35{{n}\choose{3}}+35{{n}\choose{4}}+21{{n}\choose{3}}+7{{n}\choose{2}}+n$.
\end{enumerate}
\ec 
\bpr Let $L$ be an $n$-Lie algebra of dimension $d=n$ with the basis $X=\{x_1,x_2,\dots,x_n\}$ such that $x_n>x_{n-1}>\dots>x_2>x_1$.

(1). It is clear from the definition of basic commutators of weight $1$. 
 
(2). By Theorem \ref{weight2}, we have 
\[l_d^n(2)=l_n^n(2)={{n}\choose{n}}=\dfrac{n!}{(n-n)!\times n_!}=1.\]

(3). We know that the general form of the basic commutators of weight $3$ is as follows:
\[[[x_{11},\dots,x_{1n}],x_{22},\dots,x_{2n}]\in L^3,\]
where $x_{ij}\in X$, for all $1\leq j\leq n$ and $i=1,2$. By part $(2)$, we know that at $n=d$, there is only one basic commutator of weight $2$ as $[x_n,x_{n-1},\dots,x_2,x_1]$.  Therefore, it is enough to determine only the status of the outer components of the second bracket. Since $n=d$, so $(n-1)$-components $x_{22},\dots,x_{2n}$ are equal to $(n-1)$-components of inner bracket. Without reducing the totality, assume that $x_{1j}=x_{2j}$, for all $2\leq j\leq n$. Hence in this situation, according to the  Jacobi identity and skew-symmetric property of $n$-Lie algebras, we have 
\begin{align*}
&[[x_{11},\dots,x_{1n}],x_{22},\dots,x_{2n}]\\
&=[[x_{11},\dots,x_{1n}],x_{12},\dots,x_{1n}]\\
&=[[x_{11},\dots,x_{1n}],x_{12},\dots,x_{1n}]
+[x_{11},\cancelto{0}{[x_{12},x_{12}\dots,x_{1n}]},x_{13},\dots,x_{1n}]\\
&\quad+[x_{11},x_{12},\cancelto{0}{[x_{13},x_{12}\dots,x_{1n}]},x_{14},\dots,x_{1n}]
+\dots
+[x_{11},x_{12},\dots,\cancelto{0}{[x_{1n},x_{12}\dots,x_{1n}]}]\\
&=[[x_{11},\dots,x_{1n}],x_{12},\dots,x_{1n}].
\end{align*}  
Thus, to count the number of all  basic commutators of weight $3$ on $d$ letters (which $d=n$), 
it is enough to calculate the number of all possible modes to select $(n-1)$ distinct letters from $n$ letters (in other words, it is equal to the number of all basic commutators of weight $2$ on $d$ letters in  $(n-1)$-Lie algebras), which according to Theorem \ref{weight2}, is equal to
\[l_n^d(3)=l_n^n(3)={{n}\choose{n-1}}=n.\]
So the number of all basic commutators of weight $w=3$, where $n=d$, is equal to $n$. 

(4). We know that the general form of basic commutators of weight $4$ can be considered as $[[[x_{11},x_{12},\dots,x_{1n}],x_{22},\dots,x_{2n}],x_{32},\dots,x_{3n}]\in L^4$. Note that $x_{22},\dots,x_{2n}$
 and $x_{32},\dots,x_{3n}$ can be considered as strings such as  basic commutators of weight $2$ on $n$ letters and with length $n-1$. In this case, according to the definition of the basic commutators, $(x_{12},\dots,x_{1n})\leq(x_{22},\dots,x_{2n})\leq(x_{32},\dots,x_{3n})$ should be hold. Given that at $n=d$, there is only one basic commutator of weight $2$ as $[x_1,x_2,\dots,x_n]$, it suffices to consider the possible choices for $(x_{22},\dots,x_{2n})$ and $(x_{32},\dots,x_{3n})$. Obviously, for each $2\leq j\leq n$, $x_{2j}$ and $x_{3j}$ are equal to one of the $[x_1,x_2,\dots,x_n]$ components, and $(x_{2},\dots,x_{n})$ is smaller than or equal to $(x_{22},\dots,x_{2n})$ and $(x_{32},\dots,x_{3n})$. At $n=d$, the general form of the base shifters of weight $4$ can be rewritten as follows:    
\[[[[x_{1},x_{2},\dots,x_{n}],x_{22},\dots,x_{2n}],x_{32},\dots,x_{3n}]\in L^4,\qquad x_{2j},x_{3j}\neq x_1,~~\text{for all } 2\leq j\leq n.\]
According to the above description, it can be divided into the following two categories:
\begin{align}
&[[[x_{1},x_{2},\dots,x_{n}],x_{22},\dots,x_{2n}],x_{22},\dots,x_{2n}],\quad (x_{22},\dots,x_{2n})=(x_{32},\dots,x_{3n})\label{form1}\\ 
& \Longrightarrow {\text{the number of all basic commutators in this form is equal to~}}n .\nonumber\\ \nonumber \\
&[[[x_{1},x_{2},\dots,x_{n}],x_{22},\dots,x_{2n}],x_{32},\dots,x_{3n}],\quad (x_{22},\dots,x_{2n})\neq(x_{32},\dots,x_{3n}),\label{form2}\\
& \Longrightarrow {\text{the number of all basic commutators in this form is equal to~}}{{n}\choose{2}}.\nonumber
\end{align}
Note that the number of all basic commutator of weight $2$ on $d=n$ letters in $(n-1)$-Lie algebras is equal to $n$. Therefore, the number of all basic commutators of weight $4$ on $d=n$ letters is equal to the sum of the above two cases \eqref{form1} and \eqref{form2}; that is,
\[l_n^n(4)={{l_n^{n-1}(2)}\choose{2}}+n={{{{n}\choose{n-1}}}\choose{2}}+n={{n}\choose{2}}+n.\]

(5). Similar to part (4), we have the following categories: 
\begin{align}
&[[[x_{1},x_{2},\dots,x_{n}],x_{22},\dots,x_{2n}],x_{22},\dots,x_{2n}],x_{22},\dots,x_{2n}],\nonumber \\ 
&{\text {where}}~(x_{22},\dots,x_{2n})=(x_{i2},\dots,x_{in}),\quad i=3,4,\nonumber\\ 
& \Longrightarrow {\text{the number of all basic commutators in this form is equal to~}}n. \label{form1-w=5}\\ \nonumber \\
&[[[x_{1},x_{2},\dots,x_{n}],x_{22},\dots,x_{2n}],x_{22},\dots,x_{2n}],x_{32},\dots,x_{3n}]\quad{\text {and}}\nonumber\\
&[[[x_{1},x_{2},\dots,x_{n}],x_{22},\dots,x_{2n}],x_{32},\dots,x_{3n}],x_{32},\dots,x_{3n}], \nonumber \\
&{\text{where}}~(x_{22},\dots,x_{2n})\neq(x_{32},\dots,x_{3n}), \nonumber \\
&\Longrightarrow {\text{the number of all basic commutators in this form is equal to~}}2\times{{n}\choose{2}}.\label{form2-w=5}\qquad\qquad\qquad \\ \nonumber\\ 
&[[[x_{1},x_{2},\dots,x_{n}],x_{22},\dots,x_{2n}],x_{32},\dots,x_{3n}],x_{42},\dots,x_{4n}], \nonumber \\
&{\text{where}}~(x_{i2},\dots,x_{in})\neq(x_{j2},\dots,x_{jn}),\quad i,j=2,3,4, \nonumber \\ 
&\Longrightarrow {\text{the number of all basic commutators in this form is equal to~}}{{n}\choose{3}}.\label{form3-w=5}\qquad\qquad\qquad\qquad \ 
\end{align}
Hence, $l_n^n(5)$ is equal to the sum of the three cases \eqref{form1-w=5}, \eqref{form2-w=5}, and \eqref{form3-w=5} as follows:
\begin{align*}
l_n^n(5)&={{l_n^{n-1}(2)}\choose{3}}+2{{l_n^{n-1}(2)}\choose{2}}+n\\
&={{{{n}\choose{n-1}}}\choose{2}}+2{{{n}\choose{n-1}}\choose{2}}+n\\
&={{n}\choose{3}}+2{{n}\choose{2}}+n.
\end{align*}
\indent By a similar process to the previous one, it can be proved that 
\begin{align*}
l_n^n(6)&={{{{n}\choose{n-1}}}\choose{4}}+3{{{n}\choose{n-1}}\choose{3}}+3{{{n}\choose{n-1}}\choose{2}}+n
={{n}\choose{4}}+3{{n}\choose{3}}+3{{n}\choose{2}}+n,\qquad\qquad\ 
\end{align*}
\begin{align*}
l_n^n(7)&={{{{n}\choose{n-1}}}\choose{5}}+4{{{n}\choose{n-1}}\choose{4}}+6{{{{n}\choose{n-1}}}\choose{3}}+4{{{{n}\choose{n-1}}}\choose{2}}+n\\
&={{n}\choose{5}}+4{{n}\choose{4}}+6{{n}\choose{3}}+4{{n}\choose{2}}+n,\\ \\
l_n^n(8)&={{{{n}\choose{n-1}}}\choose{6}}+5{{{{n}\choose{n-1}}}\choose{5}}+10{{{{n}\choose{n-1}}}\choose{4}}+10{{{{n}\choose{n-1}}}\choose{3}}+5{{{n}\choose{n-1}}\choose{2}}+n\\
&={{n}\choose{6}}+5{{n}\choose{5}}+10{{n}\choose{4}}+10{{n}\choose{3}}+5{{n}\choose{2}}+n,\\ \\
l_n^n(9)&={{{{n}\choose{n-1}}}\choose{7}}+6{{{{n}\choose{n-1}}}\choose{6}}+15{{{{n}\choose{n-1}}}\choose{5}}+20{{{{n}\choose{n-1}}}\choose{4}}+15{{{{n}\choose{n-1}}}\choose{3}}+6{{{n}\choose{n-1}}\choose{2}}+n\\
&={{n}\choose{7}}+6{{n}\choose{6}}+15{{n}\choose{5}}+20{{n}\choose{4}}+15{{n}\choose{3}}+6{{n}\choose{2}}+n,\\ \\
l_n^n(10)&={{{{n}\choose{n-1}}}\choose{8}}+7{{{{n}\choose{n-1}}}\choose{7}}+21{{{{n}\choose{n-1}}}\choose{6}}+35{{{{n}\choose{n-1}}}\choose{5}}+35{{{{n}\choose{n-1}}}\choose{4}}+21{{{{n}\choose{n-1}}}\choose{3}}\\
&\quad+7{{{n}\choose{n-1}}\choose{2}}+n\\
&={{n}\choose{8}}+7{{n}\choose{7}}+21{{n}\choose{6}}+35{{n}\choose{3}}+35{{n}\choose{4}}+21{{n}\choose{3}}+7{{n}\choose{2}}+n.
\end{align*}
\epr

According to Corollary \ref{Cor-l_n^n(w)}, we obtain $l_n^n(w)=\displaystyle\sum_{i=1}^{w-2}a_i{{n}\choose{i}}$, and in Table  \ref{tab3} (for $i\geq 1$),  see the order between the coefficients.

\begin{table}[h!]
\caption{coefficients ${{n}\choose{i}}$'s}\label{tab3}\begin{tabular}{|c|c|c|c|c|c|c|c|c|}\hline
 &$a_1$ &$a_2$ &$a_3$ &$a_4$ &$a_5$ &$a_6$ &$a_7$ &$a_8$\\ \hline
$l_n^n(4)$& 1 &1 &-- &-- &-- &-- &-- &-- \\ \hline  
$l_n^n(5)$& 1 &2 &1 &-- &-- &-- &-- &-- \\ \hline
$l_n^n(6)$& 1 &3 &3 &1 &-- &-- &-- &-- \\ \hline  
$l_n^n(7)$& 1 &4 &6 &4 &1 &-- &-- &-- \\ \hline
$l_n^n(8)$& 1 &5 &10 &10 &5 &1 &-- &-- \\ \hline  
$l_n^n(9)$& 1 &6 &15 &20 &15 &6 &1 &-- \\ \hline
$l_n^n(10)$& 1 &7 &21 &35 &35 &21 &7 &1 \\ \hline  
\end{tabular}
\end{table}

According to Table \ref{tab3}, we assume that the following relations are established.
\begin{itemize}
\item 
$a_1=a_{w-2}=1$, \qquad
$a_2=a_{w-3}=w-3$, \qquad and \qquad
$a_i=a_{w-i-1}$, for each $i$.
\item 
Comparing the coefficients in Table \ref{tab3} with the values of the Khayyam triangle, it is observed that to calculate the coefficients $a_i$ in $l_n^n(w)$, it is sufficient to obtain the $(w-3)$th power coefficients of the binomial, that is, $(x+y)^{w-3}$. In other words, $a_i$ is the coefficient $x^{w-3-(i-1)}y^{i-1}$. Each  $a_i$, for $1\leq i\leq w-2$, can also be considered the $i$th digit in the answer of $(11)^{w-3}$. 
\end{itemize}

Table \ref{tab4} denotes the number of basic commutators in different weights with $n=d$, which can be obtained by the  formulas in Corollary \ref{Cor-l_n^n(w)} and manual calculations.

\begin{table}[h!]
\caption{Number of basic commutators in different weights}\label{tab4}\begin{tabular}{|c|c|c|c|c|c|c|c|c|c|c|}\hline 
$n=d$&2& 3 & 4 & 5 & 6 & 7 & 8 & 9 & 10\\ \hline 
$w=1$&2& 3 & 4 & 5 & 6 & 7 & 8 & 9 & 10\\ \hline
$w=2$&1& 1 & 1 & 1 & 1 & 1 & 1 & 1 & 1 \\ \hline
$w=3$&2& 3 & 4 & 5 & 6 & 7 & 8 & 9 & 10\\ \hline
$w=4$&3& 6 & 10 & 15 & 21 & 28 & 36 & 45 & 55 \\ \hline
$w=5$&6& 10 & 20 & 35 & 56 & 84 & 120 & 165 & 220 \\ \hline
$w=6$&9& 15 & 35 & 70 & 126 & 210 & 330 & 495 & 715 \\ \hline
$w=7$&18& 21 & 56 & 126 & 252 & 462 & 792 & 1287 & 2002 \\ \hline
$w=8$&30& 28 & 84 & 210 & 462 & 924 & 1716 & 3003 & 5005 \\ \hline
$w=9$&56& 36 & 120 & 330 & 792 & 1716 & 3432 & 6435 & 11440 \\ \hline
$w=10$&99& 45 & 165 & 495 & 1287 & 3003 & 6435 & 12870 & 24310 \\ \hline
\end{tabular}\end{table}
We obtain the following simple formulas that are easily proved by induction:
\begin{itemize}
\item 
$l_n^n(4)=l_n^n(3)+(l_n^n(3)-1)+(l_n^n(3)-2)+\dots+2+1=\sum\limits_{i=0}^{l_{n}^{n}(3)-1}(l_n^n(3)-i)=\Big{(}l_{n}^{n}(3)\times l_{n}^{n}(3)\Big{)}+\sum\limits_{i=1}^{l_{n}^{n}(3)}i=\Big{(}l_{n}^{n}(3)\Big{)}^2+\displaystyle\dfrac{l_n^n(3)(l_n^n(3)-1)}{2}=\Big{(}l_{n}^{n}(3)\Big{)}^2+{{l_n^n(3)}\choose{2}},$ where $n\geq2$.
\item 
$l_n^n(w)=\sum\limits_{w'=2}^wl_{n-1}^{n-1}(w')=l_n^n(w-1)+l_{n-1}^{n-1}(w),$ where $n>3$.
\item 
$l_3^3(w)=l_3^3(w-1)+(w-1)$, where $w\geq 3$, $l_3^3(1)=3$, and $l_3^3(2)=1$.  
\item 
$l_4^4(w)=\begin{cases} 
l_4^4(w-1)+(w-1)r,\quad w=2r,\\ \\ 
l_4^4(w-1)+wr,\quad w=2r+1.
\end{cases}$.
\end{itemize}
We know that the number of all commutators (all of the basic and nonbasic commutators) with weight $w$ is $d^{m_w}$. Hence,
\[\underbrace{\stackrel{d}{\bigcirc}\stackrel{d}{\bigcirc}\dots\stackrel{d}{\bigcirc}\stackrel{d}{\bigcirc}}_{n-times}.\]
Also, it is obvious that the number of all nonbasic commutators of weight $w$, denoted by $L_d^n(w)$, is equal to
\begin{equation}\label{nonbasicweightw}
L_d^n(w)=d^{m_w}-l_d^n(w).
\end{equation}
On the other hand, by  Theorem \ref{weight2}, we have 
$l_d^n(2)={{d}\choose{n}}$, and so the number of all nonbasic commutators of weight 2 is equal to 
\begin{equation}\label{nonbasicweight2}
L_d^n(2)=d^{m_2}-l_d^n(2)=d^n-{{d}\choose{n}}.
\end{equation}
Now, we calculate the number of all basic commutators of weight 3. Every commutator of weight 3 has $n+(n-1)=2n-1$ components, in which $n$ of them form the commutators of weight 2, and others are weight 1 and form the exterior part in the bracket.    
It is obvious that every nonbasic commutator of weight 2 forms a nonbasic commutator of weight 3, independent of the selection of the $n-1$ components in external parts. So these must be removed, too. 
We denote the number of all of such cases by ${L'}_d^n(3)$, which is equal to 
\begin{equation}\label{nonbasicweight3-part1}
{L'}_d^n(3)=L_d^n(2)\times d^{n-1}.
\end{equation}
On the other hand, note that if a basic commutator of weight 2 is arranged in a nonbasic form way with the external part, then the result will be nonbasic again. So we should remove such cases. The number of these, denoted by ${L''}_d^n(3)$, is calculated as follows:
\begin{align}
{L''}_d^n(3)&=l_d^n(2)\times L_d^{n-1}(2) \label{nonbasicweight3-part2-0}\\
&={{d}\choose{n}}\times\left(d^{n-1}-{{d}\choose{n-1}}\right).\label{nonbasicweight3-part2}
\end{align} 
Thus, 
\begin{equation}\label{nonbasicweight3}
L_d^n(3)={L'}_d^n(3)+{L''}_d^n(3).
\end{equation} 
Moreover, it may be that two parts are basic, but the third condition of the definition of basic commutators is not true. Hence commutators in this situation must be replaced by terms with the basic form, and the number of them must be decreased. In this case, we can replace it by using the Jacobi identity, such that  every term of the Jacobi identity has the  basic form, or  otherwise, repeat  this replacing for nonbasic terms in the Jacobi identity. Note that if any term of Jacobi identity does not have a basic form, then we use the Jacobi identity next-handed to replace it (i.e., if the first time, we write the right Jacobi identity, then next time, must be used the left Jacobi identity to replace the other terms). 
Note that by attention to the collecting process, all of the above terms can be assumed.
We denote the number of them with $L_w^*(n,d,m_{w-1},n-1)$, where $m_w=m_{w-1}+n-1$. It is important that if $n=d$, then each product of two basic commutators as $ab$, in which $a$ is a basic commutator of weight $w-1$ and $b$ is a basic string of length $n-1$, has a basic form with weight $w$. Thus 
\begin{equation}\label{L*_w(n,n,m_{w-1},n-1)}
L_w^*(n,n,m_{w-1},n-1)=0.
\end{equation}
In general, 
\begin{equation}\label{L_w*(n,d,m_{w-1},n-1)}
L_w^*(n,d,m_{w-1},n-1)=l_d^n(w-1)\times l_d^{n-1}(2)-\kappa =l_d^n(w-1)\times{{d}\choose{n-1}}-\kappa.
\end{equation}  
We are going to find $\kappa$ by induction on $w$. \\
Let us first show that the number of all basic commutators of weight $w$ is equal to $\kappa$. 
According to the above explanation, it is clear that the number of all of the nonbasic commutators is equal to the sum of both $L^n_d(w)$ and $L_w^*(n,d,m_{w-1},n-1)$, and so other commutators are basic. Thus the number of all basic commutators of weight $w$ is 
\begin{align*}
l_d^n(w)&=d^{m_w}-L_d^n(w)-L_w^*(n,d,m_{w-1},n-1)\\
&=d^{m_w}-({L'}_d^n(w)+{L''}_d^n(w))-L_w^*(n,d,m_{w-1},n-1)\\
&=d^{m_w}-\left(L_d^n(w-1)\times d^{n-1}\right)-\left(l_d^n(w-1)\times L_d^{n-1}(2)\right)-\left(l_d^n(w-1)\times l_d^{n-1}(2)-\kappa\right)\\
&=d^{m_w}-\left((d^{m_{w-1}}-l_d^n(w-1))\times d^{n-1}\right)-\left(l_d^n(w-1)\times (d^{n-1}-l_d^{n-1}(2))\right)\\
&\quad-\left(l_d^n(w-1)\times l_d^{n-1}(2)-\kappa\right)\\
&=d^{m_w}-\underbrace{(d^{m_{w-1}}\times d^{n-1})}_{d^{m_w}}+(l_d^n(w-1)\times d^{n-1})-\left(l_d^n(w-1)\times d^{n-1}\right)\\
&\quad+(l_d^n(w-1)\times l_d^{n-1}(2))-\left(l_d^n(w-1)\times l_d^{n-1}(2)\right)+\kappa
=\kappa.
\end{align*}
Now, we calculate $\kappa$ by induction on $w$.
It is obvious that the smallest value of $w$ must be $3$, 
and so, for $w=3$, we have $l_d^n(2)={{d}\choose{n}}$, $l_d^{n-1}(2)={{d}\choose{n-1}}$, and 
\begin{equation}\label{L_3*(n,d,m_{2},n-1)}
L_3^*(n,d,m_{2},n-1)={{d}\choose{n}}{{d}\choose{n-1}}-\sum_{i=1}^{d-n+1}\ \sum_{j=i+1}^{d-1}
(d-j)\left[{{d-i+1}\choose{n-1}}-j+i+1\right].
\end{equation}  
Therefore, by \eqref{nonbasicweight2}, \eqref{nonbasicweight3-part1}, \eqref{nonbasicweight3-part2}, and \eqref{L_3*(n,d,m_{2},n-1)}, we have
\begin{align*}
l_d^n(3)=\kappa&=\sum_{i=1}^{d-n+1}\ \sum_{j=i+1}^{d-1}
(d-j)\left[{{d-i+1}\choose{n-1}}-j+i+1\right] .
\end{align*}
A comparison of two equations \eqref{basicweight3} and \eqref{L_3*(n,d,m_{2},n-1)} shows that there are no other basic commutators of weight $3$ except  equation \eqref{basicweight3}. Therefore, we obtain a formula to count basic commutators of weight $3$ in any given $n$-Lie algebras with arbitrary dimension $d\geq n$, and hence the following theorem is proved. 
\begin{theorem}\label{main-formula-weight3}
Let $F$ be a free $n$-Lie algebra over the ordered set $X=\{x_i|x_{i+1}>x_i,~i=1,2,\dots,d\}$. Then the number of basic commutators of weight $3$ is 
\begin{equation}\label{basicweight3}
l_d^n(3)=\sum_{i=1}^{d-n+1}\ \sum_{j=i+1}^{d-1}
(d-j)\left[{{d-i+1}\choose{n-1}}-j+i+1\right] .
\end{equation}
\end{theorem}
By a similar way and induction, we can rewrite \eqref{nonbasicweight3-part1}, \eqref{nonbasicweight3-part2-0}, and \eqref{nonbasicweight3} and obtain the following equations:
\begin{align*}
{L'}_d^n(w)=L_d^n(w-1)\times d^{n-1},\ \ \   
&{L''}_d^n(w)=l_d^n(w-1)\times L_d^{n-1}(2), \  
&L_d^n(w)={L'}_d^n(w)+{L''}_d^n(w). 
\end{align*}

\begin{figure}[h!]
\caption{}\label{fig1}
\includegraphics[scale=0.4]{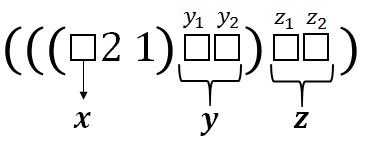}
\end{figure}

Now, in another way, we calculate the number of basic commutators of weights $4$ and $5$ on $d\geq 4$ letters. Then we generalize it to the state of arbitrary $w$ and obtain the final and general formula for the various arbitrary values of $n$, $d$, and $w$. 
\begin{theorem}\label{main-formula-weight4}
Let $F$ be a free $n$-Lie algebra over the ordered set $X=\{x_i|x_{i+1}>x_i,~i=1,2,\dots,d\}$. 
 Then the number of basic commutators of weight $4$ is 
\begin{equation}\label{basicweight4}
l_d^n(4)=\sum_{j=1}^{\alpha_0}\beta_{j^*}\left({{{{d}\choose{n-1}}}\choose{2}}+{{d}\choose{n-1}} \right),
\end{equation}
where $\alpha_0={{d-1}\choose{n-1}}$, and if ${{k-1}\choose{n-1}}+1\leq j\leq {{k}\choose{n-1}}$ (for $k=n-1,n,n+1,n+2,\dots,d-1$), then $j^*={{k-1}\choose{n-1}}+1$ and $\beta_{j^*}=(d-n-j^*+2)$.
\end{theorem}
\bpr
Let $d\geq n$, let $X=\{x_i~|~x_{i}>x_{i-1},~i=1,2,3,4,\dots,d\}$, and let $w=4$. We know that $l_{d-1}^{n-1}(2)={{d-1}\choose{n-1}}$ is the number of collected strings of length $n-1$, too. We calculate the collected strings of length $n-1$ and arrange them as follows. Then we number them from right to left (We use the $j$ counter to count them. Obviously, $j=1,2,\dots,l_{d-1}^{n-1}(2)$). Note that here we display $n=3$ just to make it easier to display the components, but the same procedure applies to the general case $n$. Hence
\begin{equation}\label{d>n,w=4-part1}
\begin{tabular}{ccccccccc}
$x_{d-1}x_{d-2}$&$x_{d-1}x_{d-3}$&$\dots$&$x_4x_3$&$x_4x_2$&$x_4x_1$&$x_3x_2$&$x_3x_1$&$x_2x_1$\\
$\downarrow$&$\downarrow$& &$\downarrow$&$\downarrow$&$\downarrow$&$\downarrow$&$\downarrow$&$\downarrow$\\
$j={{d-1}\choose{n-1}}$&$j={{d-1}\choose{n-1}}-1$&$\dots$&$j=6$&$j=5$&$j=4$&$j=3$&$j=2$&$j=1$
\end{tabular}
\end{equation}
The general form of the basic commutators of weight $4$ is as follows: 
\[(((a_1,a_2,a_3),b_2,b_3),c_2,c_3),\qquad a_r,b_s,c_q\in X.\]
First, we select the first collected string from the chain \eqref{d>n,w=4-part1} (for $j=1$) and replace $a_2,a_3$ with it; see Figure \ref{fig1}.
We have to calculate the number of all basic commutators by counting the number of allowed modes for cells such as $x$, $y$, and $z$. Since the corresponding term $j=1$ in the chain \eqref{d>n,w=4-part1} is the first and smallest sentence, the allowed states for selecting $x$ are $x_3,x_4, \ldots ,x_{d-1}$, and $x_d$, so the number of these states is equal to $d-(n-1)=d-2$. 

Now, we want to form a chain of all basic commutators of weight $2$ (with length $n-1$) on $d$ letters $x_1,x_1,\dots,x_{d}\in X$. We know that $l_{d-1}^{n-1}(2)={{d-1}\choose{n-1}}$ is the number of all basic commutators of weight $2$ and length $n-1$. We calculate the basic commutators of weight $2$ and then arrange them by the defined relation ``$<$" on basic commutators as follows. Then we label them from right to left (We use the counter $i$, to count them. Obviously, $i=1,2,\dots,l_{d}^{n-1}(2)$). Hence 
\begin{equation}\label{d>n,w=4-part2}
{\footnotesize{\underbrace{x_dx_{d-1}}_{\substack{\downarrow\\i={{d}\choose{n-1}}}}>\underbrace{x_dx_{d-2}}_{\substack{\downarrow\\i={{d}\choose{n-1}}-1}}>\underbrace{x_{d-1}x_{d-2}}_{\substack{\downarrow\\i={{d}\choose{n-1}}-2}}>\dots >\underbrace{x_dx_2}_{\substack{\downarrow\\i=2d-3}}>\dots >\underbrace{x_4x_2}_{\substack{\downarrow\\i=d+1}}>\underbrace{x_3x_2}_{\substack{\downarrow\\i=d}}>\underbrace{x_dx_1}_{\substack{\downarrow\\i=d-1}}>\dots >\underbrace{x_3x_1}_{\substack{\downarrow\\i=2}}>\underbrace{x_2x_1}_{\substack{\downarrow\\i=1}}.
}}
\end{equation}
Also, since the corresponding sentence $j=i=1$ in the chain \eqref{d>n,w=4-part2} is the smallest sentence, all the sentences in the chain \eqref{d>n,w=4-part2} can be placed after it. On the other hand, both $z\geq y$ and $y\geq x_2x_1$ conditions must be met. Obviously, both sentences we choose from the chain \eqref{d>n,w=4-part2} for $y$ and $z$, one of which is smaller than the other, which can be considered as $y$ and another as $z$. Therefore, the number of choice states $y$ and $z$, provided that $y$ and $z$ are distinct, is equal to ${{{{d}\choose{n-1}}}\choose{w-2}}={{{{d}\choose{n-1}}}\choose{2}}$. If $y=z$, then we have ${{d}\choose{n-1}}$ choices for them. 
Thus, in this case $j=1$, the number of basic commutators is equal to $(d-2)\times\left[{{{{d}\choose{n-1}}}\choose{2}}+{{d}\choose{n-1}}\right]$. 

By the same process, we also calculate the number of basic commutators in the cases $j=2,3,4,\dots,{{d-1}\choose{n-1}}$, and the sum of the number of basic commutators in the various cases $j$ is equal to the number of all the basic commutators of weight $4$ on $d$ letters. Moreover,
\[l_d^n(4)=\sum_{j=1}^{\alpha_0}\beta_{j^*}\left({{{{d}\choose{n-1}}}\choose{2}}+{{d}\choose{n-1}} \right),\]
where $\alpha_0={{d-1}\choose{n-1}}$, and if ${{k-1}\choose{n-1}}+1\leq j\leq {{k}\choose{n-1}}$, (for $k=n-1,n,n+1,n+2,\dots,d-1$), then put $j^*={{k-1}\choose{n-1}}+1$ and $\beta_{j^*}=(d-n-j^*+2)$. Therefore it is clear that $1\leq j^*\leq d-n+1$.  \epr

In the following, we want to get a formula for calculating the number of basic commutators of any weight $w$ in every $d$-dimensional $n$-Lie algebras. For this purpose, it is necessary to note that in the case of duplicate selection, the placement of strings can also have several different modes, the number of which must also be calculated. For example, suppose $w=5$. We have to count the number of allowed choices for $x$, $y$, and $z$. If all three strings $x$, $y$, and $z$ are distinct, then they are arranged as follows: 
\begin{align*}
((((a_1,\underbrace{a_2,\dots, a_n}_{a}),\underbrace{x_2,x_3,\dots, x_n}_{x}),\underbrace{y_2,y_3,\dots, y_n}_{y}),\underbrace{z_2,z_3,\dots, z_n}_{z}), \qquad z> y> x\geq a.
\end{align*} 
Indeed if two of the strings $x$, $y$, and $z$ are the same, then the following two states must be considered and counted:
\begin{align*}
&((((a_1,\underbrace{a_2,\dots, a_n}_{a}),\underbrace{x_2,x_3,\dots, x_n}_{x}),\underbrace{x_2,x_3,\dots, x_n}_{y}),\underbrace{z_2,z_3,\dots, z_n}_{z}),\qquad z>y=x\geq a,\\
&((((a_1,\underbrace{a_2,\dots, a_n}_{a}),\underbrace{x_2,x_3,\dots, x_n}_{x}),\underbrace{y_2,y_3,\dots, y_n}_{y}),\underbrace{y_2,y_3,\dots, y_n}_{z}), \qquad z= y> x\geq a.
\end{align*} 
Also, if $x=y=z$, then 
\[((((a_1,\underbrace{a_2,\dots, a_n}_{a}),\underbrace{x_2,x_3,\dots, x_n}_{x}),\underbrace{x_2,x_3,\dots, x_n}_{y}),\underbrace{x_2,x_3,\dots, x_n}_{z}),\qquad z=y=x\geq x.\]
In the weight $w=4$, there exist two outer strings $x$ and $y$, and so there were two choice modes for $x$ and $y$, which are as follows: $(1)$-A state where both $x$ and $y$ are distinct (Equivalently, $x$ and $y$ are replaced by exactly two distinct strings). $(2)$-The state that $x$ and $y$ are equal (Equivalently, $x$ and $y$ are replaced by only one string). Also, in the case $w=5$, there exist three outer strings $x$, $y$, and $z$. Hence, there were three choice modes for $x$, $y$, and $z$, which are as follows: $(1)$-A state that all three strings $x$, $y$, and $z$ are distinct  (Equivalently, $x$, $y$, and $z$ are replaced by exactly three distinct strings, in which case there is only one arrangement for these three strings as $z>y>x$). $(2)$-The state that two strings of $x$, $y$, and $z$ are equal (Equivalently, $x$, $y$, and $z$ are replaced by exactly two distinct strings, in which case there is only two arrangement for these strings as $z=y>x$ and $z>y=x$). $(3)$-The state where all of the three strings are equal together (Equivalently, $x$, $y$, and $z$ are replaced by only one string, in which case there is only one arrangement for strings as $z=y=x$).

In general, for an arbitrary $w$, the number of outer strings is equal to $w-2$, called $y_1,y_2,y_3,\dots,y_{w-2}$ (where it can replaced by basic strings of length $n-1$, or in other words, by basic commutators of weight $2$ in $(n-1)$-Lie algebras) on $d$ letters, in which the number of all them is equal to $l_{d}^{n-1}(2)={{d}\choose{n-1}}$). If we want to replace them with $w-2$ distinct strings (a state in which they are all distinct), then we have only one arrangement as $y_1<y_2<y_3<\dots<y_{w-2}$. If we replace all of them with just one string (a state in which they are all repetitive), then we have one arrangement $y_1=y_2=y_3=\dots=y_{w-2}$.   
If we replace them with $w-3$ distinct strings (a state in which exactly two strings are same) or replace them with $2$ distinct strings (a state in which some of the $y_i$ are equal to string $1$ and others are equal to string $2$), then the number of arrangements of them is equal to $\alpha_2=\alpha_{w-1}={{w-3}\choose{0}}={{w-3}\choose{w-3}}$. 

Continuing this process, we can say that if we want to replace them with $i$ distinct strings or replace them with $w-i-1$ (a state in which exactly $i$ strings are the same), then we have $\alpha_i$ different arrangements, where    $\alpha_i=\alpha_{w-i+1}={{w-3}\choose{i-2}}={{w-3}\choose{w-i+1}}$ and $2\leq i\leq w-1$.    

Now, considering the above explanations and the proof of the previous theorem, we are ready to introduce the formula for calculating the number of basic commutators. 
\begin{theorem}\label{main-formula-weightw}
Let $F$ be a free $n$-Lie algebra over the ordered set $X=\{x_i|x_{i+1}>x_i:~i=1,2,\dots,d\}$ and let $w$ be a positive integer number. Then the number of basic commutators of weight $w$ is 
\begin{equation}\label{basicweightw}
l_d^n(w)=\sum_{j=1}^{\alpha_0}\beta_{j^*}\left(\sum_{i=2}^{w-1}\alpha_i{{{{d}\choose{n-1}}}\choose{w-i}}\right),
\end{equation}
where $\alpha_0={{d-1}\choose{n-1}}$, $\alpha_i$, ($2\leq i\leq w-1$) is the coefficient of the $(i-2)$th sentence in Newton's binomial expansion $(a+b)^{w-3}$, and if ${{k-1}\choose{n-1}}+1\leq j\leq {{k}\choose{n-1}}$, (for $k=n-1,n,n+1,n+2,\dots,d-1$), then $j^*={{k-1}\choose{n-1}}+1$ and $\beta_{j^*}=(d-n-j^*+2)$.

\end{theorem}
\bpr
According to the recent explanation, it is proved by a similar process of the proof of  Theorem \ref{basicweight4}. \epr

\section{Relationship between the number of basic commutators of $n$-Lie algebras and Lie algebras}
 
We know that the $n$-Lie algebras are independent of the selection of the bracket, and the concept of basic commutators is also related to the free $n$-Lie algebras. In addition, if $L$ is a Lie algebra over the field $\Bbb F$ with the bracket $[-,-]: L \times L\longrightarrow L$, then for each natural number $n\in \Bbb N$, $L$ can be equipped with an $n$-linear function $[-,\dots,-]: \underbrace{L\times L \times\dots\times L}_{n-times}
\longrightarrow L$, such that $L$ forms an $n$-Lie algebra with the bracket $[-,\dots,-]$ on the field $\Bbb F$ in which
\begin{align*}
[-,\dots,-]: \underbrace{L\times L \times\dots\times L}_{n-times}
&\longrightarrow L\\ 
(x_1, x_2, \dots, x_n)&\longmapsto [x_1, x_2, \dots, x_n]:= [[\dots [[x_1, x_2], x_3], \dots ]].
\end{align*}
The above explanation may raise the question of whether there is a connection between $l_d(w)$ and $l^n_d(w)$, which are the number of basic commutators of weight $w$ and on $d$ generators in Lie algebras and $n$-Lie algebras?

In this section, we intend to answer this question and examine the relationship between them.
By using \eqref{Witt-formula}, we have
\begin{equation*}
l^2_d(w)=l_d(w)=\dfrac{1}{w}\sum_{s|w}\mu(s)d^{\frac{w}{s}},
\end{equation*} 
for all $w,d\in \Bbb N$.  

On the other hand, we know that
\begin{align*}
{{d}\choose{2}}&=\dfrac{d(d-1)}{2!}=\dfrac{d^2-d}{2!}
=\dfrac{1}{2}\sum_{s|2}\mu(s)d^{\frac{2}{s}}=l_d^2(2).
\end{align*} 
Now, let $n=3$. Then 
\begin{align*}
{{d}\choose{3}}&=\dfrac{d^3-d}{3!}-\dfrac{3d^2-3d}{3!}\\
&=\dfrac{d^3-d}{3\times 2!}-\dfrac{d^2-d}{2!}\\
&=\dfrac{1}{2!}\left(\dfrac{1}{3}\sum_{s|3}\mu(s)d^{\frac{3}{s}}\right)-\dfrac{1}{2}\sum_{s|2}\mu(s)d^{\frac{2}{s}}\\
&=\dfrac{1}{2!}l_d(3)-l_d(2).
\end{align*}
By continuing this process, Table \ref{tab5} can be obtained for different values of $n\in\Bbb N$. 



Now, we assume that the following equation is true for each natural number smaller
than or equal to $n_0$:
\begin{equation}\label{induction-assumpsion}
{{d}\choose{n}}=\sum_{s=2}^n(-1)^{n-s}\left(
\dfrac{a_s}{(s-1)!\times (s + 1)(s + 2)\ldots (n-1)n}\right)l_d(s),
\end{equation}
where, $a_s>0$, for all $2\leq s\leq n$. Then for $n_0+1$, we have
\begin{align*}
{{d}\choose{n_0+1}}&=\dfrac{\prod\limits_{i=0}^{d-1}(d-i)}{(d-(n_0+1))!\times (n_0+1)!} \qquad\qquad\qquad\qquad\qquad\qquad\qquad\qquad\qquad\qquad\qquad 
\end{align*}
\begin{align*}
&=\dfrac{\prod\limits_{i=0}^{n_0}(d-i)}{(n_0+1)!}\\ \\
&=\dfrac{(d-n_0)\prod\limits_{i=0}^{n_0-1}(d-i)}{(n_0+1)\times n_0!}\\ \\
&=\dfrac{(d-n_0)}{n_0+1}\times \dfrac{\prod\limits_{i=0}^{n_0-1}(d-i)}{n_0!}\\ \\
&=\dfrac{(d-n_0)}{n_0+1}\times\left(\sum_{s=2}^{n_0}(-1)^{n_0-s}\left(\dfrac{a_s}{(s-1)!\times (s+1)(s+2)\dots (n_0-1)n_0}\right)l_d(s)\right) \\ \\
&=\sum_{s=2}^{n_0}\dfrac{(-1)^{n_0-s}(d-n_0)a_s}{(s-1)!\times (s+1)(s+2)\dots (n_0-1)n_0(n_0+1)}l_d(s) \\ \\
&=\sum_{s=2}^{n_0}\dfrac{(-1)^{n_0-s}(d-n_0)a_s}{(s-1)!\times (s+1)(s+2)\dots (n_0-1)n_0(n_0+1)}\left(\dfrac{1}{s}\sum_{r|s}\mu(r)d^{s/r}\right) \\ \\
&=\sum_{s=2}^{n_0}\dfrac{(-1)^{n_0-s}(d-n_0)a_s}{(s-1)!\times (s+1)(s+2)\dots (n_0-1)n_0(n_0+1)}\left(\dfrac{1}{s}\sum_{r|s}\mu(s/r)d^{r}\right) \\  \\
&=\sum_{s=2}^{n_0}\dfrac{(-1)^{n_0-s}a_s}{(s-1)!\times (s+1)(s+2)\dots (n_0-1)n_0(n_0+1)}\left(\dfrac{(d-n_0)}{s}\sum_{r|s}\mu(s/r)d^{r}\right)   \\ \\
&=\sum_{s=2}^{n_0}\dfrac{(-1)^{n_0-s}a_s}{(s-1)!\times (s+1)(s+2)\dots (n_0-1)n_0(n_0+1)}\left(\dfrac{1}{s}\sum_{r|s}\mu(s/r)(d-n_0)d^{r}\right) \\ \\
&=\sum_{s=2}^{n_0}\dfrac{(-1)^{n_0-s}a_s}{(s-1)!\times (s+1)(s+2)\dots (n_0-1)n_0(n_0+1)}\left(\dfrac{1}{s}\sum_{r|s}\mu(s/r)d^{r+1}\right) 
\end{align*}
\begin{align*}
&\quad+\sum_{s=2}^{n_0}\dfrac{(-1)^{n_0-s}a_s}{(s-1)!\times (s+1)(s+2)\dots (n_0-1)n_0(n_0+1)}\left(\dfrac{1}{s}\sum_{r|s}\mu(s/r)(-n_0)d^{r}\right) \\  \\
&=\sum_{s=2}^{n_0}\dfrac{(-1)^{n_0-s}a_s}{(s-1)!\times (s+1)(s+2)\dots (n_0-1)n_0(n_0+1)}\left(\dfrac{1}{s}\sum_{r|s}\mu(s/r)d^{r+1}\right)  \\  \\
&\quad+\sum_{s=2}^{n_0}\dfrac{(-1)^{n_0-s}(-n_0)a_s}{(s-1)!\times (s+1)(s+2)\dots (n_0-1)n_0(n_0+1)}\left(\dfrac{1}{s}\sum_{r|s}\mu(s/r)d^{r}\right)  \\  \\
&=\sum_{s=2}^{n_0}\dfrac{(-1)^{n_0-s}a_s}{(s-1)!\times (s+1)(s+2)\dots (n_0-1)n_0(n_0+1)}\left(\dfrac{1}{s}\sum_{r|s}\mu(s/r)d^{r+1}\right)  \\  \\
&\quad-\sum_{s=2}^{n_0}\dfrac{(-1)^{n_0-s}n_0a_s}{(s-1)!\times (s+1)(s+2)\dots (n_0-1)n_0(n_0+1)}l_d(s) \\ \\
&=\sum_{s=2}^{n_0+1}\dfrac{(-1)^{n_0-s+1}b_s}{(s-1)!\times (s+1)(s+2)\dots (n_0-1)n_0(n_0+1)}l_d(s).\\
\end{align*}
Thus equation \eqref{induction-assumpsion} holds for each natural number $n$, and we have
\begin{equation}\label{connection}
{{d}\choose{n}}=\sum_{s=2}^n(-1)^{n-s}\left(
\dfrac{a_s}{(s + 1)(s + 2) \ldots (n-1)n}\right)l_d(s).
\end{equation}
Moreover,
\begin{align}
{{d}\choose{n}}&=\sum_{s=2}^n(-1)^{n-s}\left(
\dfrac{a_s}{3\times 4\times \dots\times (s-1)(s + 1)(s + 2) \ldots (n-1)n}\right)l_d(s)\nonumber\\
&=\sum_{s=2}^n(-1)^{n-s}\left(
\dfrac{a_s}{\prod\limits_{\substack{r=3\\ r\neq s}}^{n}r}\right)l_d(s).\label{true and main connection}
\end{align}

{\footnotesize{
\begin{sidewaystable}
\begin{center}\vspace{16cm}
\begin{tabular}{|c|c|c|c|c|c|c|c|c|c|c|}\hline 
 &$l_d(2)$ &$l_d(3)$ &$l_d(4)$ &$l_d(5)$ &$l_d(6)$ &$l_d(7)$ &$l_d(8)$ &$l_d(9)$ &$l_d(10)$ \\ \hline
${{d}\choose{2}}=\dfrac{1}{2!}\prod\limits_{i=0}^1(d-i)$&$\frac{1}{1}$ &$-$ &$-$ &$-$ &$-$ &$-$ &$-$ &$-$ &$-$  \\ \hline
${{d}\choose{3}}=\dfrac{1}{3!}\prod\limits_{i=0}^2(d-i)$&$\frac{-3}{3}$ &$\frac{1}{2!}$ &$-$ &$-$ &$-$ &$-$ &$-$ &$-$ &$-$ \\ \hline
${{d}\choose{4}}=\dfrac{1}{4!}\prod\limits_{i=0}^3(d-i)$&$\frac{12}{4\times 3}$ &$\frac{-6}{4\times 2!}$ &$\frac{1}{3!}$ &$-$ &$-$ &$-$ &$-$ &$-$ &$-$ \\ \hline
${{d}\choose{5}}=\dfrac{1}{5!}\prod\limits_{i=0}^4(d-i)$&$\frac{-60}{5\times 4\times 3}$ &$\frac{35}{5\times 4\times 2!}$ &$\frac{-10}{5\times 3!}$ &$\frac{1}{4!}$ &$-$ &$-$ &$-$ &$-$ &$-$ \\ \hline
${{d}\choose{6}}=\dfrac{1}{6!}\prod\limits_{i=0}^5(d-i)$&$\frac{360}{6\times 5\times 4\times 3}$ &$\frac{-224}{6\times 5\times 4\times 2!}$ &$\frac{85}{6\times 5\times 3!}$ &$\frac{-15}{6\times 4!}$ &$\frac{1}{5!}$ &$-$ &$-$ &$-$ &$-$  \\ \hline
${{d}\choose{7}}=\dfrac{1}{7!}\prod\limits_{i=0}^6(d-i)$&$\frac{-2520}{7\times 6\times\dots\times 3}$ &$\frac{1603}{7\times 6\times 5\times 4\times 2!}$ &$\frac{-735}{7\times 6\times 5\times 3!}$ &$\frac{175}{7\times 6\times 4!}$ &$\frac{-21}{7\times 5!}$ &$\frac{1}{6!}$ &$-$ &$-$ &$-$  \\ \hline
${{d}\choose{8}}=\dfrac{1}{8!}\prod\limits_{i=0}^7(d-i)$&$\frac{20160}{8\times 7\times\dots\times 3}$ &$\frac{-12810}{8\times 7\times\dots\times 4\times 2!}$ &$\frac{6770}{8\times 7\times 6\times 5\times 3!}$ &$\frac{-1960}{8\times 7\times 6\times 4!}$ &$\frac{322}{8\times 7\times 5!}$ &$\frac{-28}{8\times 6!}$ &$\frac{1}{7!}$ &$-$ &$-$ \\ \hline
${{d}\choose{9}}=\dfrac{1}{9!}\prod\limits_{i=0}^8(d-i)$&$\frac{-181440}{9\times 8\times\dots\times 3}$ &$\frac{113589}{9\times 8\times\dots\times 4\times 2!}$ &$\frac{-67320}{9\times 8\times\dots\times 5\times 3!}$ &$\frac{22449}{9\times 8\times 7\times 6\times 4!}$ &$\frac{-4536}{9\times 8\times 7\times 5!}$ &$\frac{546}{9\times 8\times 6!}$ &$\frac{-36}{9\times 7!}$ &$\frac{1}{8!}$ &$-$ \\ \hline
${{d}\choose{10}}=\dfrac{1}{10!}\prod\limits_{i=0}^9(d-i)$&$\frac{1814400}{10\times 9\times\dots\times 3}$ &$\frac{-1109472}{10\times 9\times \dots\times 4\times 2!}$ &$\frac{724550}{10\times 9\times\dots\times 5\times 3!}$ &$\frac{-269324}{10\times 9\times\dots\times 6\times 4!}$ &$\frac{63273}{10\times 9\times 8\times 7\times 5!}$ &$\frac{-9450}{10\times 9\times 8\times 6!}$ &$\frac{870}{10\times 9\times 7!}$ &$\frac{-45}{10\times 8!}$ &$\frac{1}{9!}$ \\ \hline
\end{tabular}\caption{The connection between $l_d^2(w)$ and $l^n_d(2)={{d}\choose{n}}$}\label{tab5}
\end{center}
\end{sidewaystable}
}}

\newpage
According to the above description, the following theorem is proved easily.
\bt \label{conection-weight2}
The number of basic commutators of weight $2$ on $d$ generators in $n$-Lie algebras can be concluded by the number of basic commutators of weights $w\leq n$ in Lie algebras.
\et 
\bpr
It is easy to check from equation \eqref{connection} and Theorem \ref{weight2}. \epr

The following corollary expresses the relationship between the number of basic commutators in Lie algebras and $n$-Lie algebras. The following formula is Witt's generalized formula.
\bc 
The number of basic commutators of arbitrary weight $w$ on $d$ generators in $n$-Lie algebras, can be concluded by the number of basic commutators in Lie algebra. We have 
\[l_d^n(w)=\sum_{j=1}^{\alpha_0}\beta_{j^*}\left(\sum_{i=2}^{w-1}\alpha_i\left(\sum_{s=2}^{w-i}(-1)^{w-i-s}\left(
\dfrac{a_s}{\prod\limits_{\substack{r=3\\ r\neq s}}^{w-i}r}\right)l_{d^*}(s)\right)\right),
\]
where $d^*={{d}\choose{n-1}}$, $\alpha_0={{d-1}\choose{n-1}}$, $\alpha_i$ ($2\leq i\leq w-1$) is the coefficient of the $(i-2)$th sentence in Newton's binomial expansion $(a+b)^{w-3}$. If ${{k-1}\choose{n-1}}+1\leq j\leq {{k}\choose{n-1}}$, (for $k=n-1,n,n+1,n+2,\dots,d-1$), then $j^*={{k-1}\choose{n-1}}+1$ and $\beta_{j^*}=(d-n-j^*+2)$.
\ec 
\bpr
It is proved from equation \eqref{basicweightw} in Theorem \ref{main-formula-weightw} and Theorem \ref{conection-weight2}. \epr

The following theorem is the most important and practical result of this paper.  
\bt \label{F^i/F^j}
Let $F$ be a free $n$-Lie algebra and let $F^i$ be the $i$th term of the lower central series of $F$, for each $i\in\Bbb N$. Then $\dfrac{F^i}{F^{i+c}}$ is  abelian of dimension $l_d^n(i)+l_d^n(i+1)+\dots+l_d^n(i+c-1)$, where $c=0,1,2,\ldots,i$.   
\et 
\bpr It is proved similar to its version in the theory of Lie algebras. \epr 

\ \\
{\small 
Farshid Saeedi$^*$: Corresponding author\\ 
Department of Mathematics, Mashhad Branch, Islamic Azad University, Mashhad, Iran.\\
E-mail address: saeedi@mshdiau.ac.ir \\ \ \\
Seyedeh Nafiseh Akbarossadat \\
Department of Mathematics, Mashhad Branch, Islamic Azad University, Mashhad, Iran.\\
E-mail address: n.akbarossadat@gmail.com  
}

\end{document}